\documentclass[11pt]{article}
\usepackage{a4}
\usepackage{epsfig}
\usepackage[latin3]{inputenc}
\usepackage{amsmath,amsthm,amssymb}
\oddsidemargin=\evensidemargin
\def\runninghead#1#2{\pagestyle{myheadings}
\markboth{{\protect\footnotesize\it{\quad #1}}\hfill}
{\hfill{\protect\footnotesize\it{#2\quad}}}} \headsep=15pt

\newcommand{\R}{{\rm{I\! R}}}

\newcommand{\ep}{\varepsilon}

\newtheorem{Theorem}{Theorem}[section]
\newtheorem{Corollary}[Theorem]{Corollary}

\newtheorem{Lemma}[Theorem]{Lemma}
\newtheorem{Remark}[Theorem]{Remark}
\newtheorem{Proposition}[Theorem]{Proposition}

\newcommand{\DN}{\Delta_{i} N}
\newcommand{\DNu}{\Delta_{1} N}

\newcommand{\DJu}{\Delta_{i} J_1}

\newcommand{\DW}{\Delta_{i} W}

\newcommand{\DX}{\Delta_{i} X}
\newcommand{\DXz}{\Delta_{i} X_0}
\newcommand{\DXu}{\Delta_{i} X_1}
\newcommand{\DXd}{\Delta_{i} \tilde J_2}
\newcommand{\DXq}{(\Delta_{i} X)^2}
\newcommand{\DXquarta}{(\Delta_{i} X)^4}
\newcommand{\DXzq}{(\Delta_{i} X_0)^2}
\newcommand{\DXuq}{(\Delta_{i} X_1)^2}
\newcommand{\DXdq}{(\Delta_{i} \tilde J_2)^2}
\newcommand{\DXduq}{(\Delta_{1} \tilde J_2)^2}
\newcommand{\IVXqleqr}{ I_{\{ \DXq\leq r(h)\}}}

\newcommand{\IVXuqleqqr}{ I_{\{ \DXuq\leq 4r(h)\}}}

\newcommand{\intIi}{\int_{t_{i-1}}^{t_i}}
\newcommand{\Ii}{]t_{i-1},t_i]}

\newcommand{\à}{\`a}

\newcommand{\vp}{\vspace{0.3cm}}
\newcommand{\vm}{\vspace{0.5cm}}
\newcommand{\vu}{\vspace{1cm}}

\newcommand{\dlim}{{\rm dlim}}
\newcommand{\Plim}{{\rm Plim}}
\newcommand{\ce}{\centerline}

\newcommand{\ri}{\rightarrow}

\newcommand{\dimo}{\noindent {\it Proof.  }}
\newcommand{\n}{\noindent}
\newcommand{\ba}{\begin{array}{c}}
\newcommand{\ea}{\end{array}}
\newcommand{\bteo}{\begin{Theorem} }
\newcommand{\eteo}{\end{Theorem} }
\newcommand{\bcor}{\begin{Corollary} }
\newcommand{\ecor}{\end{Corollary} }
\newcommand{\bprop}{\begin{Proposition} }
\newcommand{\eprop}{\end{Proposition} }
\newcommand{\blem}{\begin{Lemma} }
\newcommand{\elem}{\end{Lemma} }
\newcommand{\brem}{\begin{Remark} }
\newcommand{\erem}{\end{Remark} }
\newcommand{\boss}{\begin{Remark} }
\newcommand{\eoss}{\end{Remark} }
\newcommand{\beqlab}{\begin{equation} \label }
\newcommand{\beq}{\begin{equation} }
\newcommand{\eeq}{\end{equation} }

\begin{document}

\ce{\large\bf NON PARAMETRIC THRESHOLD ESTIMATION FOR MODELS}
\ce{\large\bf WITH STOCHASTIC DIFFUSION COEFFICIENT AND JUMPS}
\vp\ce{Cecilia Mancini} \ce{\footnotesize Universit\à di Firenze,
Dipartimento di matematica per le decisioni,}
\ce{\footnotesize via C.Lombroso, 6/17, 50134 Florence, Italy, tel. +39 055 4796808, fax  +39 055 4796800}
\ce{\footnotesize
cecilia.mancini@dmd.unifi.it}


 \vu \ce{\bf \sc \bf Summary}

\vm \n We consider a stochastic process driven by a diffusion and jumps.
 We devise a technique, which is based on a discrete record
 of observations, for identifying the times when jumps larger than a suitably
 defined threshold occurred. The technique allows also
 jump size estimation. We prove the consistency of a nonparametric estimator of
the integrated infinitesimal variance of the process continuous part when the jump component with infinite activity
is Lévy.
Central limit results are proved in the case where the  jump component has
 finite activity.
 Some simulations illustrate the reliability of the methodology in finite
 samples.\footnote{
{\it AMS 2000 subject classifications}. Primary: 62G05,
62G20, 
62M99; 
secondary: 7M05 
37M10.\\
The results of this paper were presented:
to the 3rd World Congress of the Bachelier Finance Society, Chicago, 21-24 July  2004;
to the workshop Quantitative Methods in Finance, Newton Institute, Cambridge, 31/1/05 - 4/2/05;
to the 9th World Congress of the Econometric Society, London, 19th - 24th August 2005,
http://www.econ.ucl.ac.uk/eswc2005/
} 

\vu\n {\it Key words}: discrete observations, non parametric
estimation, models with stochastic diffusion coefficient and jumps,
threshold, integrated infinitesimal variance of the continuous
component,
asymptotic properties.

\section{Introduction}

We consider a stochastic process $X$ starting from $x_0\in \R$ at
time $t=0$ and such that
\beq\label{eqY}
dX_t=a_tdt+\sigma_tdW_t+ dJ_t,
\; t\in]0,T], \eeq where $a$ and $\sigma$ are progressively
measurable processes,
$W$ is a standard Brownian motion and $J$ is a pure jump process. A
jump process is said to have {\it finite activity} (FA) when it
makes a.s. a finite number of jumps in each finite time interval,
otherwise it is said to have {\it infinite activity} (IA). We
provide an estimate of $\int_0^T \sigma^2_t dt$, denoted by $IV$,
given discrete observations $\{x_0, X_{t_1}, ..., X_{t_n}\}$. The
estimator is consistent both when $J$ has FA and when the IA
component of $J$ is Lévy. $IV$ stands for {\it Integral of the
infinitesimal Variance} of the continuous part of $X$; $\int_0^T
\sigma^2_t dt$ is also called {\it integrated volatility} in the
financial econometric literature. When $J$ has FA we also give an
estimate of jump times and sizes, while, when $J$ has IA we can
identify the instants when jumps are larger than a given threshold.
 These results have important applications in financial
econometrics, see the reviews in Andersen et al. (2005) and Barndorff-Nielsen
and Shephard (2006).\\
The method we propose here extends previous work in Mancini (2001)
and Mancini (2004) allowing for infinite jump activity and very mild
assumptions on $a$ and $\sigma$.

Nonparametric estimation of the diffusion coefficient $\sigma$ has
been studied, in absence of the jump component, e.g. by
Barndorff-Nielsen and Shephard (2002). For a review see Fan (2005).
However, the inclusion of jumps within financial  models seems to be
more and more necessary for practical applications (Das, 2002;
Piazzesi, 2005; Bates, 2002).
 In the literature on non parametric inference for stochastic processes
driven by diffusions plus jumps, several  approaches have been
proposed to separate the diffusion part and the jump part given
discrete observations.\\
Berman (1965) defined {\it power variation} estimators of the sum of
given powers of the jumps.
Recently these have been recovered and developed in
Barndorff-Nielsen and Shephard (2004a), Woerner (2006) and Jacod (2006).\\
Barndorff-Nielsen and Shephard (2004a, 2004b) define and use the
{\it bipower} and the {\it multipower variation} processes to
estimate $\int_0^T \sigma^p_t dt$ for given values of $p$, and in
particular they focus on $p=2$. They assume that $\sigma$ is
independent of the leading Brownian motion (in the financial
literature this is called {\it no leverage} assumption)  and that
the jump process has finite activity. In particular they  build a
test for the presence of jumps in the data generating process.
Barndorff-Nielsen et al. (2006) and Woerner (2006) show that, in
particular cases,  the consistency and central limit theorem of the
multipower variation estimators
can be  extended in the presence of infinite activity jump processes.\\
Bandi and Nguyen (2003) and Johannes (2004)
assume that $a_t\equiv a(X_t), \sigma_t\equiv \sigma(X_t)$ and that
$J$ has FA bounded jumps. They use Nadaraya Watson kernels to obtain
pointwise estimators of the functions  $a(x)$ and $\sigma(x)$ and
aggregate information about $J$. Mancini and Renò (2006) combine the
kernel and the threshold methods
to improve the estimation of the jump part and they extend the results
to the  infinite jump activity framework. \\
Our contribution to the extant literature can be summarized as
follows. First,
 in the FA case, threshold estimation is a
more effective way to identify  intervals $]t_{j-1}, t_j]$ where $J$
jumped. Second,
the threshold estimator of  $IV$ is more efficient (in the Cramer-Rao
inequality lower
bound sense) than the multipower variation estimators. Finally, the
consistency of the threshold estimator holds even under leverage and
when the observations are not equally spaced, both in the FA and in
the IA of jump cases. An alternative extension has been made in
Jacod (2006), where, in order
to obtain a central limit theorem, the diffusion coefficient dynamics
has to be specified.\\
The good performance of our estimator on finite samples of realistic
length is shown within three different simulated models.

An outline of the paper is as follows: in section 2 we introduce the
framework and the notations; in section 3 we deal with the case
where $J$ has FA: we show that by the threshold method we can
asymptotically  identify each instant of jump. As a consequence we
obtain threshold estimators of $\int_0^T \sigma_s^2 ds$ and of each
stochastic size of the occurred jumps. Using results in
Barndorff-Nielsen et al.  (2005) and in Barndorff-Nielsen and
Shephard (2006)  we show the asymptotic normality of  $\hat{IV}$,
whatever the dynamics for $\sigma$. Moreover we find the asymptotic
distribution of the estimation error of the sizes of jump under the
no leverage assumption and when the jump component is a compound
Poisson process. Section 4 is devoted to the case when the
underlying process contains an infinite activity Lévy jump part: in
a quite simple way we show that the threshold estimator of $IV$ is
still consistent, even under leverage and  when the observations are
not equally spaced. Section 5 shows the performance of the estimator
of  $IV$ in finite samples within three different simulated models.
Section 6 concludes.\\

\n {\bf Acknowledgements.} I'm sincerely grateful to Rama Cont, Jean
Jacod and Roberto Renò  for the important comments on this work. I
also want to thank PierLuigi Zezza. I thank Sergio Vessella and
Marcello Galeotti who supported this research by MIUR grant number
2002013279 and Progetto Strategico.

\runninghead{}{Nonparametric threshold estimation}
\section{The framework}

On the filtered probability space ($\Omega$, $({\cal F}_t)_{t\in[0,
T]}$, ${\cal F}$, P), let $W$ be a standard Brownian motion and $J$
be a pure jump process  given by $J_1+ \tilde J_2$, where $J_1$ has
FA and $\tilde J_2$ has IA and is Lèvy.
Let $\left( X_t \right)_{t\in [0,T]}$
   be a real process starting from
  $x_0\in \R$ and such that
 \beq\label{eqY}
 dX_t=a_tdt+\sigma_tdW_t+ dJ_t,
\; t\in]0,T], \eeq where $a$, $\sigma$ are progressively measurable
processes which guarantee that (\ref{eqY}) has a unique strong
solution on $[0,T]$ which is adapted and right continuous with left
limits (se e.g. Ikeda and Watanabe, 1981; Protter, 1990).
Suppose that on the finite and fixed time horizon $[0,T]$ we dispose
of a discrete record
 $\{x_0, X_{t_1},...,X_{t_{n-1}},X_{t_n}\}$ of $n+1$ observations of a realization of $X$, with
 $t_i=ih$, for a given lag $h$, $T=nh$.

When $J$ is a pure jump L\'evy process, we can always  decompose it
as the sum of the jumps larger than one and the sum of the
compensated jumps smaller than one, as follows \beq \ba
J= J_1 + \tilde J_2,\\
\\
J_{1 s}\doteq\int_0^s\int_{|x|> 1} x \mu(dt,dx ), \quad \tilde J_{2
s}\doteq  \int_0^s\int_{|x|\leq 1} x(\mu(dt,dx)-\nu(dx)dt), \ea \eeq
where $\mu$ is the Poisson random measure of the jumps of $J$,
$\tilde\mu (dt,dx)= \mu(dt,dx)-\nu(dx)dt$ is the compensated
measure, and $\nu$ is the L\'evy measure of $J$  (see Sato, 1999 or
Ikeda and Watanabe, 1981). $\tilde J_2$ is a square integrable
martingale with infinite activity of jump.
 For each $s$, $Var(\tilde J_2$$_{s})= s\int_{|x|\leq 1} x^2
\nu(dx)\doteq s\sigma^2(1)<\infty$.
 $J_1$ is a  compound Poisson process with finite activity of jump, and we can also write
$J_1$$_{s}=\sum_{i=1}^{N_s} \gamma_i$, where $N$ is a Poisson
process with constant intensity $\lambda$, jumping at times denoted
by $\left(\tau_i\right)_{\!i=1..N_T} $, and each $ \gamma_i$, also
denoted $\gamma_{\tau_i}$, is the size of the jump occurred at
$\tau_i$. The random variables $\gamma_i$ are i.i.d. and
independent of $N$.\\
More generally a FA jump process is of the form
$J_1$$_{s}=\sum_{i=1}^{N_s} \gamma_i$, where $N$ is a non explosive
counting process and the random variables $\gamma_i$ are not
necessarily i.i.d., nor independent of
$N$.\\
Denote by $\tau^{(i)}$ the first instant a jump occurs within
$]t_{i-1}, t_i]$, if $\DN\geq 1$;  by $\gamma^{(i)}$ the size 
of this first jump within $]t_{i-1}, t_i]$, if $\DN\geq 1$;  by
$\underline{\gamma}\doteq   \min_{j=1..N_T }|\gamma_j|$ .\\

\n Next section deals with the FA case where $\tilde J_2\equiv 0$,
while in section 4 we allow $J$ to have infinite
activity, where $\tilde J_2$ is Lévy.\\

\n {\bf Further notations.}\\
 For any semimartingale $Z$, let us denote by $\Delta_{i} Z$ the increment $Z_{t_i}-Z_{t_{i-1}}$ and
 by $\Delta Z_t$ the size $Z_t - Z_{t-}$ of the jump (eventually) occurred at time $t$.  \\
$[Z]$ is the quadratic variation process associated to $Z$.\\
$[Z^{(h)}]_T$ is the estimator $\sum_{i=1}^n (\Delta_{i} Z)^2$ of the quadratic variation $[Z]_T$.\\
${\cal F}^Z$ denotes the sigma-algebra generated by the process $Z$.\\
$H.W$ is the process given by the stochastic integral $\int_0^{\cdot} H_s dW_s$.\\
$IV_t=\int_0^t \sigma^2_u du$.\\
$IQ_t=\int_0^t \sigma^4_u du$. This quantity is  called in the econometric literature {\it integrated quarticity} of  $X$.\\
By $c$ (low case) we denote generically a constant.\\
$\Plim$ means "limit in probability"; $\dlim$ means "limit in distribution".\\
If $\eta$ is a r.v., $M{\cal N}(0,\eta)$ indicates the mixed
Gaussian law having characteristic function $\phi(\theta)=
 E[e^{-\frac 1 2 \eta^2 \theta^2 }]$.\\

\section{Finite activity jumps}

\subsection{Consistency}
An important variable related to $X$ and containing $IV_T$ is the
quadratic variation at $T$ \beq\label{quadVar} [X]_T=\int_0^T
\sigma^2_t dt + \int_0^T \int_{\R} x^2 \mu(dx, dt). \eeq An estimate
of  $[X]_T$ is given by $ \sum_{j=1,...,m} (X_{t_j}-X_{t_{j-1}})^2$,
since
 $\Plim_{|\pi^{(T)}|\rightarrow 0}$ $ \sum_{j=1,...,m}$ $
(X_{t_j}-X_{t_{j-1}})^2\doteq [X]_T, $ where $\pi^{(T)}$ is a finite
partition $ \{t_0=0,t_1,..., t_m=T\}$ of $[0,T]$, and
$|\pi^{(T)}|=\max_j |t_j-t_{j-1}|$.\\
We consider in this section the case in which $J$ has FA,
so that
(\ref{quadVar}) becomes
$$[X]_T=\int_0^T \sigma^2_t dt + \sum_{j=1}^{N_T}\gamma^2_{\tau_j},$$
and the quadratic variation gives us only an aggregate information
regarding both $IV$ and  the jump sizes. In order to estimate the
contribution of $\int_0^T \sigma^2_t dt$ to $[X]_T$, the key point
is to exclude the time intervals $\Ii$ where $J$ jumped. The
following theorem provides an instrument to asymptotically
identifying such intervals.

\begin{Theorem}\label{DNugualDXleqr}
Suppose that $J=\sum_{j=1}^{N_t}\gamma_j$ is a finite activity jump
process where $N$ is a non explosive counting process and the random
variables $\gamma_j$ satisfy $\forall t\in[0,T]$ $P\{\Delta
N_t\neq 0, \gamma_{N_t}=0\}=0$.
Suppose also that\\
 1) a.s.
$\limsup\limits_{h\ri 0} \frac{\sup_{i} |\intIi a_s ds|}{\sqrt{h\log \frac{1}{h}}} \leq C(\omega)<\infty$\\
2) a.s. 
$\int_0^T \sigma^2_s ds <\infty$ and
$\limsup\limits_{h\ri 0} \frac{\sup_{i} |\intIi \sigma^2_s ds|}{h} \leq M(\omega)<\infty$;\\
3) $r(h)$ is a deterministic function of the lag $h$ between the
observations, s.t.
$$\lim\limits_{h\ri 0}r(h)= 0, \mbox{ and } \lim\limits_{h\ri 0}\frac{h \log \frac{1}{h}}{ r(h)}= 0,$$
 then for {\rm P}-almost all
  $\omega$ $\exists \bar{h}(\omega)$ s.t. $\forall h\leq \bar{h}(\omega)$
  we have
\beqlab{DN=DX} \forall i=1,...,n, \;\;\;
  I_{\{ \DN=0\}}(\omega) =  \IVXqleqr(\omega).
\eeq
\end{Theorem}

\vp Assumption 3) indicates how to choose the {\it threshold}
$r(h)$. The absolute value of  the increments of any path of the
Brownian motion (and thus of a stochastic integral with respect to
the Brownian motion) tends a.s. to zero as the deterministic
function $\sqrt{2 h \ln \frac{1}{h}}$. Therefore, for small $h$,
when we find that the squared increment $\DXq$ is larger than
$r(h)>2 h \ln \frac{1}{h}$  some jumps had to be occurred.

For the proof we need the following preliminary remarks.

$\bullet$ The Paul L\'evy law for the modulus of continuity of
Brownian motion's paths (see e.g. Karatzas and Shreve, 1999, theorem
9.25) implies that
  $$
  \mbox{a.s.} \quad  \lim\limits_{h\ri 0}\sup_{i\in \{ 1,...,n \} }\frac{|\DW|}{\sqrt{2h \log\frac{1}{h}} }
  \leq 1. $$

$\bullet$
The stochastic integral $\sigma.W$ is  a time changed Brownian
motion (Revuz and Yor 2001, theorems 1.9 and 1.10): defined the
pseudo-inverse of $\left(IV_t\right)_t$, $\xi_t=\inf\{v: IV_v
>t \}$, then
 \beq\label{magiorazDsigmaW} \Delta_{i}\left(\sigma.W\right)= B_{IV_{t_{i}}}-B_{IV_{t_{i-1}}},\eeq
  where $B$ is a Brownian motion.

$\bullet$ As a consequence,  under assumptions 1) and 2) of theorem
\ref{DNugualDXleqr}, by Karatzas and Shreve (1999, theorem 9.25) and
the monotonicity of the function $x\ln \frac 1 x$ it follows that
a.s. for small $h$
$$\sup_i \frac{|\intIi
a_sds+\intIi\sigma_s dW_s|}{\sqrt{2h \log\frac{1}{h} }}\leq
\Lambda(\omega),$$
 where $\Lambda(\omega)= C(\omega)+
\sqrt{M(\omega)} + 1$ is a finite r.v..

\vspace{0.7cm} \n {\it Proof of the theorem}.  First we show that
    a.s., for small $h$, it holds that  $\forall i, \; I_{\{\DN =0\}}\leq
  I_{\{\DXq \leq r(h)\}}$, then we will see that a.s., for small $h$, it holds also that
 $\forall i, \; I_{\{\DN =0\}}\geq  I_{\{\DXq \leq r(h)\}}$, and that will conclude our proof.

  1) For each $\omega$ set $J_{0,h}=\{i \in \{ 1,...,n \} : \DN =0\}$: to show
  that a.s., for small $h$, $I_{\{\DN =0\}}\leq
  I_{\{\DXq \leq r(h)\}}$ it is sufficient to prove that a.s., for small $h$, $\sup_{J_{0,h}} \DXq\leq r(h).$
  To evaluate the $\sup_{J_{0,h}} \DXq$, remark that a.s.
 $$
  \sup_{i\in J_{0,h}} \frac{\DXq}{r(h)} =
  \sup_{ J_{0,h}}
    \left(\frac{|\intIi a_sds+\intIi\sigma_s dW_s|}{\sqrt{2h \log\frac{1}{h} }}\right)^2 \cdot
  \frac{2h \log\frac{1}{h}}{r(h)}\leq
\Lambda^2\frac{2h \log\frac{1}{h}}{r(h)}\ri 0
  $$
 In particular,
    for small $h$, $ \sup_{i\in J_{0,h}} \frac{\DXq}{r(h)}\leq 1$,
as we need.

  2) Now we establish the other inequality. For any $\omega$
  set $J_{1,h}=\{i\in \{ 1,...,n \}: \DN \neq 0\}$. In order to prove that
 a.s., for small $h$,
 $\forall i, \; I_{\{\DN =0\}}\geq  I_{\{\DXq \leq r(h)\}}$ it is sufficient to show that a.s.,
 for small $h$,   $\inf_{i \in J_{1,h}} \DXq > r(h)$.
   In order to evaluate $\inf_{i \in J_{1,h}} \frac{\DXq}{r(h)}$ remark that
  $$\forall i\in J_{1,h},
  \frac{\DXq}{r(h)}=\frac{\left(\intIi a_s ds+\Delta_{i}
  \sigma.W\right)^2}{2h \log\frac{1}{h}}\frac{2h \log\frac{1}{h}}{r(h)}+$$
  $$+2\frac{\intIi a_s ds+\Delta_{i}
  \sigma.W}{\sqrt{r(h)}}\frac{\sum_{\ell=1}^{\DN} \gamma_{\ell}}{\sqrt{r(h)}}+
  \frac{(\sum_{\ell=1}^{\DN} \gamma_{\ell})^2}{r(h)}$$
  the first term tends a.s. to zero uniformly with respect to $i$. Since
for small $h$ we have that $\DN\leq 1$ for
  each $i$, then the other terms become
$$\frac{\gamma_{\tau^{(i)}}}{\sqrt{r(h)}}\left[ \frac{\intIi a_s ds+\Delta_{i}
  (\sigma.W)}{\sqrt{r(h)}}+
 \frac{\gamma_{\tau^{(i)}}}{\sqrt{r(h)}}\right].$$
 The contribution of the first term within brackets tends
a.s. to zero uniformly on $i$. Note that the assumption on $J$
guarantees that $P\{\underline{\gamma}=0\}=0$,
 thus a.s.
$$\lim_h \inf_{i \in J_{1,h}} \frac{\DXq}{r(h)}= \lim_h \frac{\gamma^2_{\tau^{(i)}}}{r(h)}\geq
\lim_h \frac{\underline{\gamma}^2}{r(h)}=+\infty.$$

\vspace{-0.8cm}\qed

\vspace{0.8cm} \n {\bf Remarks}.

i) Assumption  1) simply asks for the sequence $\sup_{i} |\intIi a_s
ds|/\left(h\log \frac{1}{h}\right)^{1/2}$
 keeping bounded as $h\ri 0$. It is satisfied if, for example, $\left(a_s(\omega)\right)_s$
is bounded  pathwise on $[0,T]$. In particular assumption 1) is
satisfied in a model with mean reverting drift $a_s = k\theta
-kX_s$.

If we assume that in equation (\ref{eqY}) $a$ and $\sigma$ are
processes having right continuous paths with left limits (c\àdl\àg),
then assumptions 1) and 2) are immediately satisfied, since a.s.
such paths are bounded on [0,T].

 ii) Note that a FA  L\évy process satisfies that $P\{\underline{\gamma}=0\}=0$, since
$\nu\{0\}=0$ (jumps occurring with zero size are not jumps). E.g.
this is the case for a compound Poisson process with Gaussian sizes
of jump.

iii) Frequently, in practice, the lag $\Delta t_i\doteq t_i-t_{i-1}
$ between the observations of an available record $\{x_{0},
X_{t_1},...,X_{t_{n-1}},X_{t_n}\}$ is not constant  (not equally
spaced observations). Theorem \ref{DNugualDXleqr}, and thus also
(\ref{consistIntVol}) below, is still valid. In fact if we set
$h\doteq \max_i \Delta t_i$,  all the fundamental ingredients of the
proof of theorem \ref{DNugualDXleqr} hold:
$$\lim\limits_{h\ri 0}\sup_{i\in \{ 1,...,n \} }\frac{|\DW|}{\sqrt{2h \log\frac{1}{h}} }\leq
 \lim\limits_{h\ri 0}\sup_{i\in \{ 1,...,n \} }\frac{|\DW|}{\sqrt{2\Delta t_i \log\frac{1}{\Delta t_i}}} \leq 1, $$
 by  the monotonicity of $x \ln \frac{1}{x}.$ Moreover, using (\ref{magiorazDsigmaW}),
it still holds that $$ \mbox{a.s.} \quad \sup_{i\in \{
1,...,n\}}\frac{|\Delta_{i}\sigma.W|}{\sqrt{2h \log\frac{1}{h}}
}\leq M(\omega),$$ since
$ \mbox{a.s.}\ \forall i \quad \Delta_{i} IV < \Delta t_i\cdot M(\omega)\leq h M(\omega).$\\
It is asymptotically equivalent to directly compare each $\DXq$ with
the relative $r(\Delta t_i)$: a.s. for small $h$ we have, for each
$i=1...,n$,
$$ I_{\{\DXq\leq r(\Delta t_i)\}} = I_{\{\DN =0\}}.$$

\vspace{-1cm}\qed

\vspace{0.7cm} Define $$\hat{IV}= \sum_{i=1}^n \DXq \IVXqleqr. $$
The consistency of $\hat{IV}$ is a consequence of theorem 3.1, which
is needed in order to asymptotically identify and exclude each jump
instant.

\bcor\label{StimIV} Under the assumptions of theorem
\ref{DNugualDXleqr} we have \beq\label{consistIntVol} {\rm
P}\!\lim\limits_{\!\!\!\!\! h\ri 0} \sum_{i=1}^n \DXq \IVXqleqr =
\int_0^T \sigma_t^2 dt.\eeq \ecor

\dimo Since a.s. for small $h$ we have $  I_{\{ \DN=0\}}=
\IVXqleqr,$ uniformly on $i$,
  then  $${\rm P}\!\lim\limits_{\!\!\!\!\! h\ri 0} \sum_i \DXq\IVXqleqr =
  {\rm P}\!\lim\limits_{\!\!\!\!\! h\ri 0} \sum_i \DXq I_{\{\DN=0\}}=$$
$$  \Plim \sum_{i=1}^n (\intIi a_sds+\intIi\sigma_s dW_s)^2
 - \Plim \sum_{i=1}^n  (\intIi a_sds+\intIi\sigma_s dW_s)^2 I_{\{\DN\neq 0\}}$$
 which coincides with $\int_0^T \sigma^2_u du,$ since
 $ \sum_{i=1}^n  (\intIi a_sds+\intIi\sigma_s dW_s)^2 I_{\{\DN\neq 0\}}\leq$
$ N_T \cdot $ $\sup_i(\intIi a_sds+\intIi\sigma_s dW_s)^2 \ri 0.$
  \qed

\subsection{Central limit theorems}

 As a corollary of theorem 2.2 in Barndorff-Nielsen et al. (2005), from our theorem
\ref{DNugualDXleqr} we obtain a threshold estimator of $\int_0^T
\sigma_t^4 dt$, which is alternative  to the power  variation
estimator. An estimate of $\int_0^T \sigma_t^4 dt$ is needed in
order to give the asymptotic law of the approximation  error
$\sum_{i=1}^n \DXq \IVXqleqr - \int_0^T \sigma_t^2$. We reach a
central limit result for $\hat{IV}$
whatever the dynamics for $\sigma$.\\

\bteo [{\it power variation} estimator: theorem 2.2 in
Barndorff-Nielsen and al. (2005), case  $r=4, s=0$]
 If $dX= a_sds + \sigma_s dW_s$, where $a$
is predictable and locally bounded and $\sigma$ is c\àdl\àg,
   then for $h\ri 0$
$$ \frac{1}{3}\frac{\sum_i \DXquarta }{h} \stackrel{P}\ri \int_0^T \sigma_t^4 dt.$$

\vspace{-1cm}\qed \eteo

\vspace{0.5cm} In the light of this result we now state the
following asymptotic properties of the threshold estimator of $IQ$.

\bprop Under the same assumptions as in theorem \ref{DNugualDXleqr}
and the assumptions of theorem 2.2 in Barndorff-Nielsen et al.
(2005)
we have that
$$ \frac{1}{3}\frac{\sum_i \DXquarta \IVXqleqr }{h} \stackrel{P}\ri  \int_0^T \sigma_t^4 dt.$$
\eprop

\dimo By theorem \ref{DNugualDXleqr}
$$\Plim \  \frac{1}{3}\frac{\sum_i \DXquarta \IVXqleqr }{h} =
\Plim \ \frac{1}{3}\frac{\sum_i \DXquarta I_{\{\DN=0\}} }{h}.$$ The
latter coincides with \beq\label{SumDX0quartadivisoh} \Plim \
\frac{1}{3}\frac{\sum_{i=1}^n (\intIi a_sds+\intIi\sigma_s
dW_s)^4}{h}, \eeq since \beq\label{SumDX0quartaIVNneq0divisoh}\Plim\
\frac{1}{3}\frac{\sum_{i=1}^n  (\intIi a_sds+\intIi\sigma_s dW_s)^4
I_{\{\DN\neq 0\}} }{h}\leq
 \Plim\   \Lambda^4 N_T \frac{(h\ln\frac{1}{h})^2}{3h}=0.\eeq
 Finally (\ref{SumDX0quartadivisoh})
coincides with $ \int_0^T \sigma_t^4 dt$, as Barndorff-Nielsen et
al. (2005) have shown. \qed


\vspace{0.5cm} Finally, as a corollary of theorem 1 in
Barndorff-Nielsen and Shephard (2006) we have the following result
of asymptotic normality
for our estimator $\hat{IV}_T.$\\

\bteo [Theorem 1 in Barndorff-Nielsen and
Shephard (2006)] 
If $dX= a_sds + \sigma_s dW_s$, where $a$ and $\sigma$ are
c\àdl\àg processes 
then, as $h\ri 0$
$$ \frac{[X^{(h)}]_T- [X]_T}{\sqrt h} \stackrel{d}\ri \sqrt{2} \int_0^T \sigma^2_u dB_u,$$
where $B$ is a Brownian motion independent of $X$ (recall from the
notations  that $[X^{(h)}]_T= \sum_{i=1}^n \DXq $).\qed \eteo

\bprop Under the assumptions of theorem \ref{DNugualDXleqr} and if
 $a$ is cadlag and locally bounded, $\sigma$ is cadlag and
    $ {\cal F}^{X}$-measurable,
then we have
$$ \frac{\sum_{i=1}^n \DXq \IVXqleqr - \int_0^T \sigma_t^2 dt}{
\sqrt{\frac{2}{3}\sum_i \DXquarta \IVXqleqr}}  \stackrel{d}\ri {\cal
N}\left(0, 1 \right).$$ \eprop

\dimo Denoting by $X_0$ the continuous process given by $X_{0 t}
=\int_0^t a_sds+\int_0^t\sigma_s dW_s $, for all $t\in [0,T]$,
similarly as in (\ref{SumDX0quartaIVNneq0divisoh})
$$\dlim\  \frac{\sum_{i=1}^n \DXq \IVXqleqr - \int_0^T
\sigma_t^2 dt}{ \sqrt{\frac{2}{3}\sum_i \DXquarta \IVXqleqr}} $$
coincides with
$$\dlim \frac{\sum_{i=1}^n  (\Delta_{i} X_0)^2 -\int_0^T \sigma_t^2 dt}
{\sqrt{\frac{2}{3}\sum_i (\Delta_{i} X_0)^4 - \frac{2}{3}\sum_i
(\Delta_{i} X_0)^4I_{\DN\neq 0} }}=$$
$$\dlim \frac{[X_0^{(h)}]_T-[X_0]_T}{\sqrt
h \sqrt{      \frac{2}{3}\frac{ \sum_i (\Delta_{i} X_0)^4 }{h} } }=
\dlim \frac{[X^{(h)}_0]_T-[X_0]_T}{\sqrt h} \frac{1}{\sqrt{2\int_0^T
\sigma_t^4 dt}}\ .$$
The first factor tends in law to $\sqrt 2 \int_0^T \sigma^2_u dB_u$,
by Barndorff-Nielsen and Shephard result (2006, theorem 1). However
(Jacod and Protter, 1998) $B$ is independent on the
whole $X$. 
Now the assumption $ \sigma \in {\cal F}^{X}$ ensures that $\sigma$
is independent of $B$ and thus, conditionally on $ \sigma$, $B$ is
again a Brownian motion and $\sqrt 2\int_0^T \sigma^2_u dB_u$ is
Gaussian with law ${\cal N}(0, 2 \int_0^T\sigma^4_u du)$. Thus,
conditionally on $ \sigma$, we have that $
\frac{[X^{(h)}_0]_T-[X_0]_T}{\sqrt h} \frac{1}{\sqrt{2\int_0^T
\sigma_t^4 dt}} \stackrel{d}\ri {\cal N}(0,1)$.
However the convergence in distribution holds even without
conditioning. \qed

\vspace{0.5cm} \n {\bf Remarks.}

i)  A comparison with the bipower variation (BPV) estimator shows
that the advantages of the non parametric threshold method are al
least two.

The threshold estimator of $IV$ is efficient (in the Cramer-Rao
inequality lower bound sense), in fact we showed  that
$\frac{\hat{IV}_T - IV_T}{\sqrt h \sqrt{\hat{IQ}_T}}$ tends in
distribution to ${\cal N}(0,2)$, while  Barndorff-Nielsen and
Shephard (2004b, p.29) show that (under the further assumption that
$X$ is a diffusion) the limit law is ${\cal N}(0, \frac{\pi^2}{4}+
\pi -3)$.
 In particular the threshold estimator is efficient (see A\"{i}t-Sahalia, 2004 for constant $\sigma$).

Moreover, since we asymptotically identify each  jump instant, we
can apply known estimation methods for diffusion processes also to
jump-diffusion processes as soon as we have eliminated the jumps
(see e.g. Mancini and Renò, 2006).

ii)  In Mancini and Renò (2006) we show that it is possible to
consider also a time varying threshold, which is
particularly important for the practical application of the
estimator.\qed

\vspace{0.4cm} By theorem \ref{DNugualDXleqr} and by the fact that,
for small $h$, the probability of more than one jump over an
interval $\Ii$ is low, it is clear that an estimator of each jump
instant
is obtained through
$$\hat \Delta_{i} N \doteq I_{\{ \DXq > r(h)\}}.$$
Moreover  a natural estimate of each realized
 jump size is given by
$$ \hat\gamma^{(i)}\doteq \DX I_{\{\DXq >r(h)\}},$$
since when a jump occurs then the contribution of $\intIi a_u du +
\intIi \sigma_u dW_u$ to $\DX$ is asymptotically negligible. In
Mancini (2004)  we have shown the consistency of each
$\hat\gamma^{(i)}$ when $T\ri \infty$. However we only gave a lower
bound for the speed of convergence when the  coefficients $\sigma$
and $a$ are stochastic processes. Here we show that, at least under
the no leverage assumption and when $J_1$ is Lévy,
the speed is exactly $\sqrt n$.\\

\bteo\label{CLTperSizesOfJump}
If   $J$ is a compound Poisson process, if a.s.
$\limsup\limits_{h\ri 0} \frac{\sup_{i} |\intIi a_s ds|}{h^\mu} \leq
C(\omega)<\infty$ for some $\mu > 0.5$ (which is the case if
$a$ is c\àdl\àg), if $\sigma$ is 
an adapted stochastic process with continuous paths,  if $\sigma$ is
independent
of $W$ and $N$,
with $E[\int_0^T \sigma^2_sds ]<\infty$, if the threshold $r(h)$ is
chosen as in theorem \ref{DNugualDXleqr} then
$$ \sqrt{n}\sum_{i} \left(\hat\gamma^{(i)}-\gamma^{(i)}I_{\{\DN \geq 1\}} \right)\stackrel{d}\ri
M{\cal N}\left(0, T\int_0^T \sigma^2_s dN_s\right). $$
\eteo

\dimo $$\sqrt{n}\sum_{i} \left(\hat\gamma^{(i)}-\gamma^{(i)}I_{\{\DN
\geq 1\}} \right)=$$
$$\sqrt{n}\sum_{i} \DX I_{\{ \DXq > r(h),\  \DN=0\}} +
\sqrt{n}\sum_{i}\! \left(\!\DX I_{\{ \DXq > r(h),\
\DN=1\}}-\gamma^{(i)}I_{\{\DN =1\}} \!\right)
$$
$$+ \sqrt{n}\sum_{i} \left(\DX I_{\{ \DXq > r(h),\  \DN\geq 2\}}-\gamma^{(i)}I_{\{\DN \geq 2\}} \right):$$
by theorem \ref{DNugualDXleqr}, a.s. for small $h$, the first term
vanishes. The third term tends to zero in probability, since
$$P\{ \sqrt{n}\sum_{i} \left(\DX I_{\{ \DXq > r(h)\}}-\gamma^{(i)}\right)I_{\{\DN \geq 2\}} \neq 0\}\leq $$
$$ P(\cup_i \{\DN \geq 2\}) \leq n O(h^2) = O(h). $$
Therefore we only have to compute the 
$$\dlim
 \sqrt{n}\sum_{i} \left(\DX I_{\{ \DXq > r(h), \DN=1\}}-\gamma^{(i)}I_{\{\DN =1\}} \right)=$$
$$\dlim
 \sqrt{n}\sum_{i} \left(\DX-\gamma^{(i)}\right)I_{\{\DN =1\}} =$$
\beq\label{pasIntermAsNormGammai}
 \dlim   \sqrt{n}\sum_{i} \left(\intIi a_s ds +\intIi \sigma_s dW_s\right) I_{\{ \DN=1\}}.\eeq
However, since for small $h$
$$P\left\{  \sqrt{n}\sum_{i} \left|\intIi a_s ds \right| I_{\{  \DN=1\}} >\ep  \right\}\leq $$
$$P\left\{  \sqrt{n}C(\omega) h^\mu N_T >\ep  \right\}\leq
P\left\{  h^{\mu-0.5}\sqrt{T}C(\omega)  N_T >\ep  \right\}\ri 0, $$
as $h\ri 0$, (\ref{pasIntermAsNormGammai}) coincides with
$$ \dlim   \sqrt{n}\sum_{i} \intIi \sigma_s dW_s I_{\{ \DN=1\}}.$$
Let us compute the characteristic function
$$ E\left[e^{i \theta \sqrt{n}\sum_{i} \intIi \sigma_s dW_s\  I_{\{ \DN=1\}} }\right]=
E\left[e^{i \theta \sqrt{T}\sum_{i} \frac{\intIi \sigma_s
dW_s}{\sqrt h}\  I_{\{ \DN=1\}} }\right]:$$ conditionally on
$\sigma$, $\frac{\intIi \sigma_s dW_s}{\sqrt h}$ are independent
Gaussian random  variables with law ${\cal N}\left(0,
\frac{\intIi\sigma^2_s ds}{h}\right)$. Since $W$ and $N$ are
independent (Ikeda and Watanabe, 1981), our characteristic function
equals
$$\Pi_{i=1}^n E\left[e^{i \theta \sqrt{T} \frac{\intIi \sigma_s dW_s}{\sqrt h}} I_{\{ \DN=1\}}
    +I_{\{ \DN\neq 1\}}  \right]=$$
\beqlab{ProdFzCar} \Pi_{i=1}^n \left( 1
+\left(e^{-\frac{1}{2}\theta^2 T\frac{\intIi\sigma^2_s ds}{h}
}-1\right)
                   e^{-\lambda h} \lambda h \right)\doteq \Pi_{i=1}^n (1+\theta_{ni}).
\eeq However $\max_i |\theta_{ni}| \ri 0$, $\sum_{i=1}^n
|\theta_{ni}| \leq \lambda$ and
$$\sum_{i=1}^n \theta_{ni}=
\lambda e^{-\lambda h}  \sum_i \left(e^{-\frac{1}{2}\theta^2
T\sigma^2_{\xi_i} }-1\right)h\ri
 \lambda\int_0^T \left(e^{-\frac{\theta^2 T}{2} \sigma^2_s} -1\right) ds,$$
 where, for each $i$, $\xi_i$ are suitable points belonging to $]t_{i-1}, t_i[$. Therefore (Chung, 1974, p.199)
(\ref{ProdFzCar}) tends to $e^{\lambda\int_0^T
\left(e^{-\frac{\theta^2 T}{2} \sigma^2_s} -1\right) ds},$ which
coincides with (Cont and Tankov, p.78) $E  \left[e^{-\frac{\theta^2
T}{2} \int_0^T \sigma^2_s dN_s}\right],$ the characteristic function
of a mixed Gaussian r.v. $\eta Z$ where $Z(P)={\cal N}(0,1)$ and
$\eta^2 = T\int_0^T \sigma^2_s dN_s.$ \qed

\section{Infinite activity jumps}

Let us now consider the case when $J$ has possibly infinite
activity. Denote
\begin{equation}\label{decompX}
\begin{array}{c}
X_{0 s}\doteq \int_0^s a_tdt+\int_0^s\sigma_tdW_t,\\
X_1\doteq X_0+J_1,\ea \eeq

\n and note that since
 $\int_{|x|\leq 1} x^2 \nu(dx)<+\infty$, then as $\ep \ri 0$ $$\sigma^2(\ep)\doteq \int_{|x|\leq\ep} x^2
\nu(dx)\ri 0.$$

\n In fact our threshold estimator is still able to extract $IV$
from the observed data. The reason is that now $$[X]_T= \int_0^T
\sigma^2_u du + \int_0^T\int_{|x|>0} x^2 \mu(dx, du)=
 \int_0^T \sigma^2_u du + \sum_{s\leq T} (\Delta J_{1 s})^2 + \sum_{s\leq T} (\Delta \tilde J_{2 s})^2 $$
 and the threshold $r(h)$ cuts off all the jumps of $J_1$ and the jumps of $\tilde J_2$ larger, in absolute value, than
 $2\sqrt{r(h)}$. However such jumps are all jumps of $\tilde J_2$ when $r(h)\ri 0$.\\

 \bteo\label{TeoStimIntVolInfActiv}
Let the assumptions 1) (pathwise boundedness condition on $a$), 2)
(pathwise boundedness  condition on $\sigma$) and 3) (choice
of the function $r(h)$) 
of theorem \ref{DNugualDXleqr} hold.
Let $J=J_1+\tilde J_2$ be such that $J_1$ has FA with $P\{\DN\neq
0\}=O(h)$ for all $h$, for all $i=1..n$; let $\tilde J_2$ be L\'evy
and be independent of $N$. Then \beq\label{goalTeo} {\rm
P}\!\lim\limits_{\!\!\!\! h\ri 0}\sum_{i=1}^n \DXq \IVXqleqr =
\int_0^T \sigma_t^2 dt. \eeq \eteo

\vspace{0.4cm}
 Jacod (2006) proves  the consistency of the threshold estimator when  the
jump process is a more general pure jump semimartingale, with the
choice $r(h)= h^{\beta}$. The proof we present here is simpler and
it allows to understand the contribution of the different jump terms
to the estimation bias. Most importantly, the advantage of the
approach presented here is that it allows to prove a central limit
theorem for $\hat{IV}$ without  any substantial assumption on
$\sigma$, while in Jacod (2006)  an assumption on the dynamics of
$\sigma$ is needed in order to get a CLT. This topic is further
developed in Cont and Mancini (2005).

To prove theorem \ref{TeoStimIntVolInfActiv} we decompose $X$ into
the sum of a jump diffusion process, $X_1$, with stochastic
diffusion coefficient and a finite activity jump part, plus an
infinite activity  compensated process $\tilde J_2$ of small jumps.
We use corollary \ref{StimIV} for the first term, and we show that
the contribution of each $\DXd$ is negligible within the truncated
version
$\sum_{i=1}^n \DXq \IVXqleqr $ of $[X^{(h)}]_T$.\\

\dimo
Since $X=X_1+\tilde J_2$, we can write
$$\left|\sum_i \DXq \IVXqleqr -\int_0^t \sigma^2 dt\right|\leq$$
\beq\label{4addendi}\begin{array}{c}
\left|\sum_i \DXuq \IVXuqleqqr -\int_0^t \sigma^2dt\right|+\\
\\
+\left|\sum_i \DXuq (\IVXqleqr - \IVXuqleqqr)\right|+\\
\\
+ 2 \left|\sum_i \DXu\DXd \IVXqleqr\right|+ \left|\sum_i \DXdq
\IVXqleqr\right|.
\end{array}\eeq

\vm \n By  corollary \ref{StimIV} we know that the first term of the
left hand side tends to zero in probability. We now
 show that the $\Plim$ of each one of the other three  terms  of the left hand side is zero.\\

Let us deal with the \underline{second term}:
$$\left|\sum_i \DXuq (\IVXqleqr - \IVXuqleqqr)\right|=$$
\beq\label{appog2} \left|\sum_i \DXuq\!\! \left( I_{\{\DXq\leq r(h),
\DXuq> 4r(h)\}}
  - I_{\{\DXq> r(h), \DXuq\leq 4r(h)\}}\right)\right|
\eeq If $I_{\{\DXq\leq r(h), \DXuq> 4r(h)\}}=1$, since
$$2\sqrt{r(h)}- |\DXd|<|\DXu|-|\DXd|\leq
|\DXu+\DXd|\leq \sqrt{r(h)},$$ then $|\DXd|>\sqrt{r(h)}$. Thus a.s.
$$\sum_i \DXuq I_{\{\DXq\leq r(h), \DXuq> 4r(h)\}}\leq \sum_i \DXuq I_{\{\DXdq > r(h)\}}\leq $$

\beq\label{appog1}  2\sum_i \left(\DXz\right)^2 I_{\{\DXdq >
r(h)\}}+ 2\sum_i \left(\sum_{j=1}^{\DN} \gamma_j\right)^2 I_{\{\DXdq
> r(h)\}}.\eeq

\n The first term is a.s.  dominated by \beq\label{appog3}
 2 \Lambda^2 h\ln \frac 1 h \sum_i I_{\{\DXdq > r(h)\}}\stackrel{P}\ri 0,
 \eeq
as $h\ri 0$, since
$$h\ln\frac{1}{h} \ n  P\{\DXdq > r(h)\}\leq h\ln\frac{1}{h}\  n\frac{E[\DXdq]}{r(h)} =
nh\sigma^2(1)\frac{h\ln\frac{1}{h}}{r(h)}\ri 0.$$
Moreover
\beq\label{PUDNneq0eDXdqger}P\left\{\sum_i (\sum_{j=1}^{\DN}
\gamma_j)^2 I_{\{\DXdq > r(h)\}} \neq 0 \right\}\leq
 P\left( \cup_i \{  \DN\neq 0, \DXdq > r(h)  \}   \right)  \leq  \eeq
$$n P\{\DNu\neq 0\}\frac{E[\DXduq]}{r(h)}=  n O(h)  \frac{h\sigma^2(1)}{r(h)}\ri 0.$$
\n For the second term of (\ref{appog2}) we note that
 by theorem \ref{DNugualDXleqr} for small $h$ on
$\{\DXuq\leq 4r(h)\}$ we have,
uniformly with respect to $i$, $\DN=0$. Therefore for small $h$
 $$\left\{\DXq> r(h), \DXuq\leq 4r(h)\right\}\subset  \left\{( \DXz+\DXd)^2> r(h)\right\}\subset $$
$$\left\{( \DXz)^2> \frac{r(h)}{4}\right\}\cup
\left\{\DXdq> \frac{r(h)}{4}\right\}. $$  However, by theorem
\ref{DNugualDXleqr}, for small $h$ a.s. $\forall i=1..n\ \   I_{\{(
\DXz)^2>\frac{r(h)}{4}\}}= 0$  and thus
$$\sum_i \DXuq I_{\{\DXq> r(h), \DXuq\leq 4 r(h)\}}\leq
\sum_i \DXzq I_{\{ \DXdq> \frac{r(h)}{4}\}}\ri 0,$$
as before in (\ref{appog3}).
Therefore (\ref{appog2}) vanishes.

\vm Let us now deal with (half) the $\Plim$ of the \underline{third
term} on the left hand side of (\ref{4addendi}), which coincides
with
\beq\label{thirdTermDoveDXdqleqr}\Plim \sum_i \DXu\DXd I_{\{
|\DX|\leq \sqrt{r(h)}, |\DXd|\leq 2\sqrt{r(h)}\}}\ .\eeq
In fact 
if  $ |\DX|\leq \sqrt{r(h)}$ and $|\DXd|>2\sqrt{r(h)}$ then
$2\sqrt{r(h)}- |\DXu|< |\DXd|-|\DXu|\leq |\DX| \leq \sqrt{r(h)}$,
i.e. $|\DXu| > \sqrt{r(h)}$, so that $|\DJu| + |\DXz|> |\DJu + \DXz|
>\sqrt{r(h)}$, \n and then \beq\label{DJ1MagRadr} \mbox{ either }
|\DJu|>\frac{\sqrt{r(h)}}{2} \mbox{ or } |\DXz| >
\frac{\sqrt{r(h)}}{2}. \eeq Since for small $h$, uniformly in $i$,
$I_{\left\{ |\DXz| > \frac{\sqrt{r(h)}}{2} \right\}}=0$, then
\beq\label{thirdTermDoveDXdqleqreDXdqgeqr} P\left\{ \sum_i
|\DXu\DXd| I_{\{ |\DX|\leq \sqrt{r(h)},\  |\DXd|>2\sqrt{r(h)}\}}
\neq 0\right\}\leq \eeq
$$P\left(\cup_i \{ |\DXd|>2\sqrt{r(h)},\ \DN\neq 0  \}\right),$$
which tends to zero as in (\ref{PUDNneq0eDXdqger}).

In order to deal now with (\ref{thirdTermDoveDXdqleqr}), note that
if $|\DX|\leq\sqrt{r(h)}$ and $|\DXd|\leq 2\sqrt{r(h)}$ then
 $$|\DJu|-|\DXz+\DXd|<|\DX|\leq\sqrt{r(h)},$$
 so a.s. $ \DN<|\DJu|<\sqrt{r(h)}+
\Lambda \sqrt{h\ln \frac{1}{h}}+|\DXd|= O\left(\sqrt{r(h)}\right), $
uniformly in $i$.
Therefore a.s. for small $h$, $\forall \ i$ on $\{ |\DX|\leq
\sqrt{r(h)}, |\DXd|$ \-$ \leq 2\sqrt{r(h)}\}$  we have
 $\DN=0$.
Thus (\ref{thirdTermDoveDXdqleqr}) is dominated by
$$ \Plim \sum_i |\DXz \DXd| I_{\{ |\DXd|\leq 2\sqrt{r(h)}\}}\leq $$
$$ \Plim \sqrt{\sum_i \DXzq} \sqrt{\sum_i \DXdq I_{\{ |\DXd|\leq 2\sqrt{r(h)}\}}}=
 \sqrt{IV_T}\ \Plim  \sqrt{\sum_i (Y^{(h)})^2}=0$$
by the Schwartz inequality and remark \ref{SuDJdsqleqrSaltiqleq9r}
below.

\vp Finally let us show that  \underline{last term} of the left hand
side of (\ref{4addendi}) tends to zero in probability. Analogously
as in (\ref{thirdTermDoveDXdqleqreDXdqgeqr})
$$ P\left\{ \sum_i \DXdq I_{\{ |\DX|\leq \sqrt{r(h)},\  |\DXd| > 2\sqrt{r(h)}\}} \neq 0\right\}\leq$$
$$P( \cup_i \{ \DN \neq 0,\  |\DXd| > 2\sqrt{r(h)}\} )\ri 0,$$
so that last term of (\ref{4addendi}) coincides with
$$ \Plim \sum_i \DXdq I_{\{ |\DX|\leq \sqrt{r(h)},\  |\DXd|\leq 2\sqrt{r(h)}\}}\leq
 \Plim \sum_i \DXdq I_{\{\DXdq\leq 4 r(h)\}}=0$$
by remark \ref{SuDJdsqleqrSaltiqleq9r}.
 \qed

\vspace{0.5cm} \brem\label{SuDJdsqleqrSaltiqleq9r} A.s., for small
$h$, uniformly in $i$, on $\{\DXdq\leq 4 r(h)\}$ we have that all
the jumps $|\Delta \tilde J_{2,s}|$ are bounded by $2 \sqrt{r(h)}$,
that is
$$(\Delta \tilde J_{2,s})^2\leq 4r(h), \quad \forall s\in ]t_{i-1}, t_i].$$
More precisely on $\{\DXdq\leq 4r(h)\}$ we have $\int_0^t \int_{ 2
\sqrt{r(h)}< |x| \leq 1}\  x \mu(ds,dx)=0$, therefore \\
$\DXd I_{\{\DXdq\leq 4 r(h)\}}$  are the increments of the process
$Y^{(h)}$ given by
$$Y^{(h)}_t\doteq \int_0^t \int_{|x|\leq 2 \sqrt{r(h)}}\  x [\mu(ds,dx)- \nu(dx)ds]
-t \int_{2\sqrt{r(h)}\leq |x|\leq 1}\ x \nu(dx).
$$
As a consequence
$$\Plim \sum_i \DXdq I_{\{\DXdq\leq 4 r(h)\}}=\Plim \sum_i(\Delta_{i} Y^{(h)})^2=$$
\beq\label{sumDXdqIVXdqleqr} \Plim\  [Y^{(h)}]_T= \Plim \int_0^T
\int _{|x|<1\wedge 2\sqrt{r(h)}}\ x^2  \mu(ds,dx)=0,\eeq \n since
last term has expectation $\sigma^2(1\wedge 2\sqrt{r(h)}\ )\ri 0$,
as $h\ri 0$.
\erem

In fact$$\frac{E[\sup_i \sum_{s\in \Ii} (\Delta \tilde J_{2,s})^2]}{r(h)}
= \frac{h}{r(h)} \sigma^2(1)\ri 0$$ as $h\ri 0$, meaning that a.s.,
uniformly on $i$, $\sum_{s\in \Ii} (\Delta \tilde J_{2,s})^2$ tends
to zero more quickly than $r(h)$, that is for small $h$, for all
$i=1..n$, $\sum_{s\in \Ii} (\Delta \tilde J_{2,s})^2<r(h),$
and thus  $(\Delta \tilde J_{2,s})^2 \leq 4 r(h)$ for each $s\in
\Ii$. \qed

\vspace{0.5cm} \n {\bf Remarks.}

i)   Everything is
still valid if we have  non equally spaced observations. 
In fact if we set, as in the remark of the previous section,
$h\doteq \max_i \Delta t_i$, the term $E[ I_{\{\DXdq > r(h)\}}]$, we
often encounter from equation (\ref{appog3}) on, is still
negligible, since as $h= \max \Delta t_i \ri 0$, $P\{\DXdq >
r(h)\}\leq \frac{E[\DXdq]}{r(h)}= \frac{\Delta t_i
\sigma^2(1)}{r(h)}\leq c\  \frac{h}{r(h)}\ $. On $\{\DXdq\leq 4
r(h)\}$ each jump size
$(\Delta \tilde J_{2,s})^2\leq 4 r(\Delta t_i)\leq 4 r(h)$, so that 
(\ref{sumDXdqIVXdqleqr}) still holds.

ii) Consistently with the results in Barnodrff-Nielsen et al. (2006)
for the multipower variations and in Jacod (2006), the asymptotic
normality of our estimator of $\int_0^T \sigma^2_s ds$ does not hold
in general if
 $X$ has an infinite activity L\'evy jump component and general cadlag coefficient $\sigma$ (Cont and Mancini 2005).
 Namely the asymptotic normality holds when the jump component has a moderate jump activity
 (when the Blumenthal-Gatoor index $\alpha$ of $J$ belongs to
 $[0,1[$), while the speed of convergence of $\hat{IV}_T$ is less than $\sqrt h$  if the  activity of jump of $J$ is
 too wild ($\alpha\in [1,2[$).

\section{Simulations}

In this section we study the performance of our threshold estimator
on finite samples. We implement the threshold estimator within three
different simulated models which are commonly used in finance: a
jump diffusion process with jump part given by a compound Poisson
process with Gaussian jump sizes; a similar model with
 stochastic diffusion coefficient correlated with the Brownian motion
driving the dynamics of $X$; and a model with an
infinite activity (finite variation) Variance Gamma jump part.   \\

MODEL 1. Let us begin with the case of a jump diffusion process with
finite activity compound Poisson jump part. We generated $N=5000$
trajectories of a process of kind
$$
dX_t= \sigma dW_t + \sum_{i=1}^{N_T} Z_i$$ with $Z_i$ i.i.d. with
law ${\cal N}(0, \eta^2)$, where $\eta=0.6$, $\sigma=0.3$ and
$\lambda$ is intentionally chosen higher than a realistic situation,
$\lambda=5$, like as in A\"{i}t-Sahalia (2004).
To generate each path we discretized EDS (\ref{eqY}) and we took
$n=6000$ equally spaced observations $X_{t_i}$ with lag $h=
\frac{1}{n}$ so that $T=1$.
We chose  $r(h)=h^{0.9}$. Figure 1 shows the distribution of the
5000 values assumed by the normalized bias term \beq\label{normBias}
\frac{\sum_{i=1}^n \DXq \IVXqleqr - \int_0^T \sigma_t^2 dt}{
\sqrt{\frac{2}{3}\sum_i \DXquarta \IVXqleqr}}\eeq versus the
standard Gaussian density (continuous line).

\begin{figure}[h]
   \begin{center}\begin{minipage}[c]{.47\linewidth}
      \begin{center}
        \includegraphics[width=4.5cm,keepaspectratio]{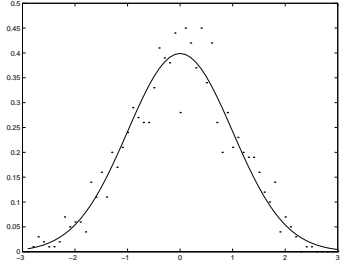}
       \end{center}
     \end{minipage}
     \caption{Distribution of 5000 values assumed by the normalized bias term (24),
     when $X$ has constant volatility and compound Poisson jumps (MODEL 1). The
continuous line is the density of the theoretical limit law ${\cal
N}(0,1)$.}
  \end{center}
\end{figure}

MODEL 2. Let us now consider  a process with jump part given by a
finite activity compound Poisson process and a stochastic diffusion
coefficient correlated with the Brownian motion driving $X$.
We generated $N=5000$ trajectories of a process of kind
$$
X_t= \ln (S_t)$$ where
$$  \frac{dS_t}{S_{t-}}= \mu dt + \sigma_t dW_t^{(1)}+ dJ_t\quad J_t= \sum_{i=1}^{N_t} Z_i,\quad
Z_i \sim {\cal N}(m_G, \nu^2),$$
$$\sigma_t =
e^{H_t},\quad d H_t= -k (H_t - \bar H) dt + \eta dW_t^{(2)},\quad
d<W^{(1)}, W^{(2)}>_t = \rho dt.
$$
Note that $$ dX_t  = (\mu - \sigma^2/2) dt + \sigma dW_t^{(1)} +
\ln(1+ \Delta J_t).$$ We chose  $\mu=0$, $\lambda=4$ and a negative
correlation coefficient $\rho= -0.7$; then we took $H_0\equiv
\ln(0.3)$, $k=1$,  $\bar H=\ln(0.25)$, $\eta=0.01$ so to ensure that
a path of $\sigma$ within $[0, T]$ varies most  between 0.2 and 0.4.
Moreover $ m_G=0.001$, $\nu=\sqrt{0.02}$ give relative amplitudes of
the jumps of $S$ most between 0.01 and 0.20. Finally we again took
$n=6000$ equally spaced observations $X_{t_i}$ with lag $h=
\frac{1}{n}$ and $r(h)=h^{0.9}$. Figure 2 shows the  distribution of
the normalized bias term (\ref{normBias}) against the asymptotic
density (continuous line).

\begin{figure}[h]
   \begin{center}\begin{minipage}[c]{.47\linewidth}
      \begin{center}
        \includegraphics[width=4.5cm,keepaspectratio]{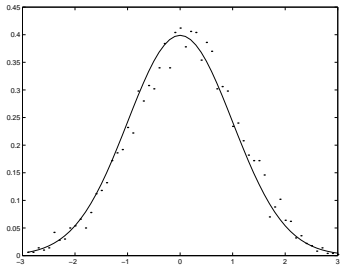}
                       \end{center}
     \end{minipage}
     \caption{Distribution of 5000 values assumed by the normalized bias term (24)
 when $X$ has stochastic volatility, negatively correlated
  with $W^{(1)}$, plus compound Poisson jumps (MODEL 2). The continuous line is the
density of the theoretical limit law ${\cal N}( 0,1)$.}
  \end{center}
\end{figure}

 MODEL 3. \n Figure 3 shows the distribution obtained in the case of
a Variance Gamma (VG) jump component. The VG process is a pure jump
process with infinite activity and finite variation. We add to it a
diffusion component $\sigma B_t$:
$$ X_t= \sigma B_t+ c G_t + \eta W_{G_t}.$$
The subordinator $G$ is a Gamma process having $Var(G_1)=b$, $B$ and
$W$ are independent Brownian motions;
we chose $N=5000$, $n=6000$ and $h=1/n$. 
$b=0.23$, $ c=-0.2$ and $\eta=0.2$ are chosen  as in Madan (2001);
$\sigma=0.3$ is chosen so that $Var(X_1)=\eta^2+c^2
b+\sigma^2=.0892$ matches the $Var (X_1)$ we obtained
for model 1. Finally $r(h)=h^{0.99}$. 

\begin{figure}[h]
\begin{center}
   \begin{minipage}[c]{.47\linewidth}
      \begin{center}
        \includegraphics[width=4.5cm, keepaspectratio]{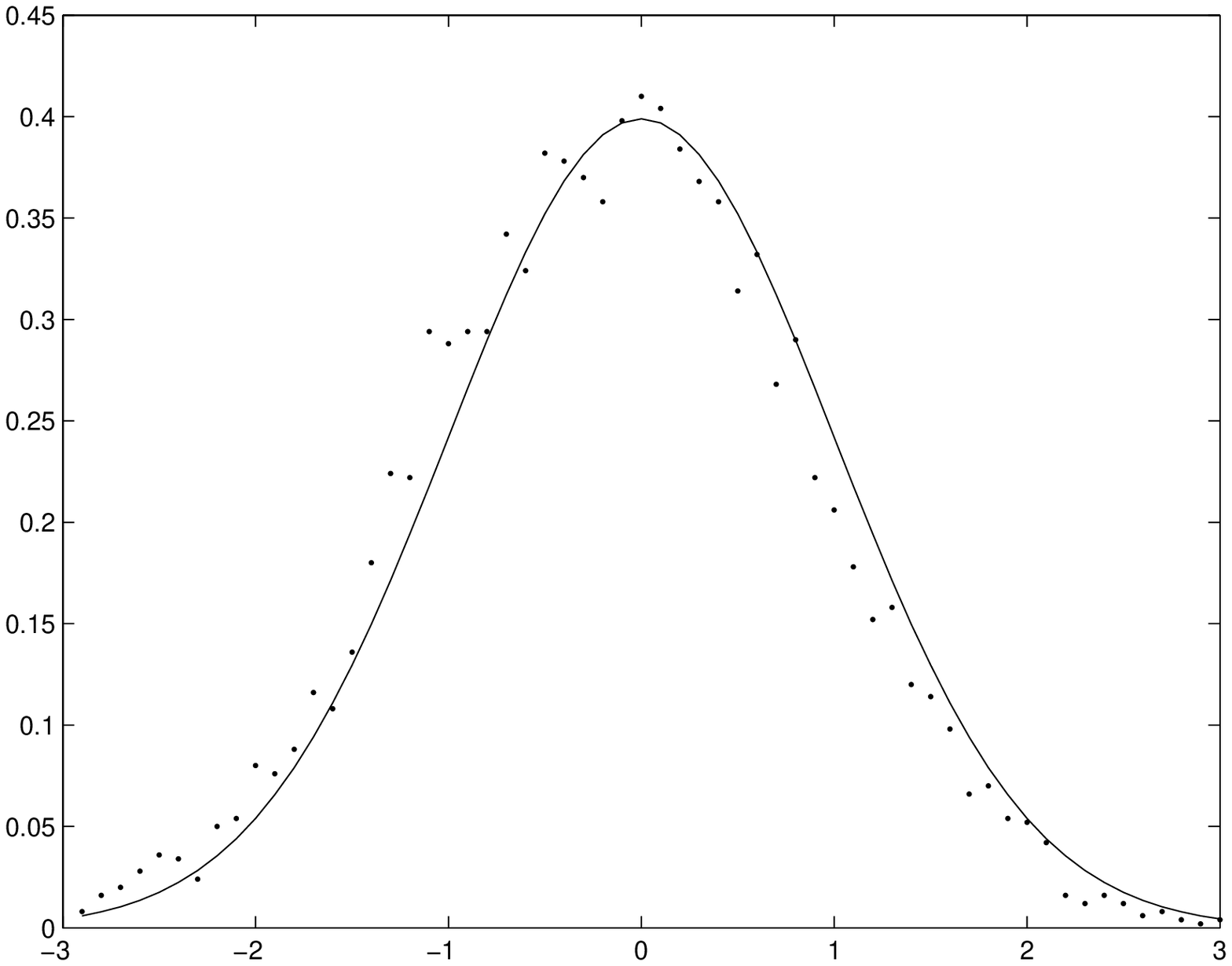}
       \end{center}
     \end{minipage}
     \caption{Distribution of 5000 values assumed by the normalized bias term (24)
      when $X$ has constant volati\-lity  plus a Variance Gamma jump part (MODEL 3). The
continuous line is the density of the theoretical limit law
${\cal N}(0,1)$.}
\end{center}
\end{figure}

\section{Conclusions}
In this paper we devise a technique for identifying the time
instants of significant jumps for a process driven by diffusion and
jumps, based on a discrete
record of observations, making use of a suitably defined threshold.\\
    We provide a consistent  estimate of $IV=\int_0^T \sigma^2_t dt$, extending previous results
    (Mancini 2001, 2004) with very mild assumptions on $a$ and $\sigma$ and, most importantly,
allowing for infinite jump activity.\\
    When $J$ has finite activity, we  give
a nonparametric estimate of the jump times and sizes, while when $J$
has a pure jump L\'evy component with infinite activity we can
identify the instants when jumps are larger than the threshold.
    When $J$ has FA we also  prove central limit results for $\hat{IV}$ and for the jump sizes estimates.

    Compared with power variations, multipower variations or kernel estimators the threshold method in the FA
case is a more effective way to identify  each interval $]t_{j-1},
t_j]$ where $J$ jumped.

We also prove that the threshold estimator of  $IV$ is efficient.\\
    Moreover, our method allows the extension of kernel estimators in diffusion frameworks to
processes driven by diffusions and jumps, provided we eliminate the jumps (Mancini and Renò, 2006).\\
    The consistency of the threshold estimator holds even under leverage, both in FA and IA cases.\\
    The threshold technique holds even when the observations are not equally spaced and also when the threshold
is time varying, which is
particularly important for the practical application of the estimator.\\

   The advantage of the approach presented here is that it allows to prove a central
limit theorem for $\hat{IV}$ without any substantial assumptions on
$\sigma$, while in Jacod (2006) an assumption on the dynamics of
$\sigma$ is needed in order to get a CLT. This topic is further
developed in Cont and Mancini (2005).

    The good performance of our estimator on finite samples of realistic length
is shown within three different simulated models.

\vu {\small \ce{\bf References}

\vm {\sc A\"{i}t-Sahalia}, Y. (2004). Disentangling volatility from
jumps. {\em Journal of Financial Economics}, 74, 487-528

{\sc Andersen}, T.G., {\sc Bollerslev}, T., {\sc Diebold}, F.X.
(2005). Parametriuc and nonparametric volatility measurement. In:
{\it Handbook of financial econometrics}, Y. A\"{i}t-Sahalia and
L.P. Hansen Eds

{\sc Bandi}, F.M., and {\sc Nguyen}, T.H. (2003). On the functional
estimation of jump-diffusion models. {\it Journal of Econometrics},
{\bf 116}, 1, pp. 293-328(36)

{\sc Barndorff-Nielsen}, O.E., {\sc Gravensen}, S.E., {\sc Jacod},
J., {\sc Podolskij}, M. and {\sc Shephard}, N. (2005). A central
limit theorem for realised power and bipower variation of continuous
semimartingales. To appear in {\em From Stochastic Analysis to
Mathematical Finance}, Festschrift for Albert Shiryaev

{\sc Barndorff-Nielsen}, O.E., {\sc Shephard} (2002). Econometric
analysis of realized volatility and its use in estimating stochastic
volatility models. {\it Journal of the Royal Statistical Society,
Series B}, 64, 253-280

{\sc Barndorff-Nielsen}, O.E. and {\sc Shephard}, N. (2004a). Power
and bipower
variation with stochastic volatility and jumps
(with discussion). {\it Journal of Financial Econometrics} {\bf 2},
1-48

{\sc Barndorff-Nielsen}, O.E. and {\sc Shephard}, N. (2004b).
Econometrics of testing for jumps in financial economics
 using bipower variation. {\it Journal of Financial Econometrics}, 2006, 4, 1-30

{\sc Barndorff-Nielsen}, O.E. and {\sc Shephard}, N.  (2006):
Variation, jumps and high frequency data in financial econometrics.
In {\em Advanced in Economics and Econometrics. Theory and
Applications}, Ninth World Congress Eds Richard Blundell, Persson
Torsten, Whitney K Newey, Econometric Society Monographs, Cambridge
University Press

{\sc Barndorff-Nielsen}, O.E., {\sc Shephard}, N.,  {\sc Winkel}, M.
(2006), Limit theorems for multipower variation in the presence of
jumps. {\it Stochastic Processes and Their Applications}, 2006, 116,
796-806

{\sc Berman}, S.M. (1965). Sign-invariant random variables and
stochastic processes with sign invariant increments. {\em Trans.
Amer. Math. Soc}, 119, 216-243

{\sc Chung}, K.L. (1974). {\it A course in probability theory}.
Academic Press Inc.

{\sc Cont}, R. and {\sc Mancini}, C. (2005).  Detecting the presence
of a diffusion and the nature of the jumps in asset prices. Working
paper

{\sc Cont}, R. and {\sc Tankov}, P. (2004). {\it Financial modelling
with jump processes}. Chapman\& Hall - CRC

{\sc Das}, S. (2002). The surprise element: jumps in interest rates.
{\it Journal of Econometrics}, 106, 27-65

{\sc Fan}, J. (2005). A selective overview of nonparametric methods
in finance. {\it Statistical Science} 20 (4), 317-337

{\sc Ikeda}, N., {\sc Watanabe}, S. (1981). {\it Stochastic
differential equations and diffusion processes}. North Holland

{\sc Jacod}, J. (2006). Asymptotic properties of realized power
variations and associated functionals of semimartingales.
arXiv, 20 April 2006 n. 0023146

{\sc Jacod}, J., {\sc Protter}, P. (1998).  Asymptotic error
distributions for the Euler method for stochastic differential
equations. {\it The Annals of Probability} {\bf 26}, 267-307

{\sc Johannes}, M. (2004). The statistical and economic role of
jumps in continuous-time interest rate models. {\it The Journal of
finance}, 59, 227-260.

{\sc Karatzas}, I., {\sc Shreve}, S.E. (1999): {\em Brownian motion
and stochastic calculus}. Springer

{\sc Madan}, D.B. (2001) Purely discontinuous asset price processes.
Advances in {\it  Mathematical Finance} Eds. J. Cvitanic,
 E. Jouini and M. Musiela, Cambridge University Press

{\sc Mancini}, C. (2001). Disentangling the jumps of the diffusion
in a geometric jumping Brownian motion. {\it Giornale dell'Istituto
Italiano degli Attuari}, Volume LXIV, Roma, 19-47

{\sc Mancini}, C., (2004).  Estimation of the parameters of jump of
a general Poisson-diffusion model. {\it Scandinavian Actuarial
Journal}, 2004, 1:42-52

{\sc Mancini}, C., {\sc Renò}, R. (2006). Threshold estimation of
jump-diffusion models and interest rate modeling. working paper

{\sc Piazzesi}, M. (2005). Bond Yields and the Federal Reserve. {\it
Journal of Political Economy} 113 (2), 311-344

{\sc Protter}, P., (1990) {\it Stochastic integration and
differential equations}. Springer-Verlag

{\sc Revuz}, D., {\sc Yor}, M. (2001) {\it Continuous martingales
and Brownian motion}. Springer

{\sc Sato}, K., (1999). {\it L\'evy Processes and infinitely
divisible distributions}. Cambridge University Press

{\sc Woerner}, J.  (2006):  Power and Multipower variation:
inference for high frequency data. In {\em Stochastic Finance}, eds
A.N. Shiryaev, M. do Ros\'ario Grossinho, P. Oliviera, M. Esquivel,
Springer, 343-364. }

\end{document}